\documentclass[12pt]{amsart}

\usepackage{amsmath}
\usepackage{amssymb}  
\usepackage{latexsym} 
\usepackage{comment}
\usepackage{url}
\usepackage{fullpage,url,amssymb,amsmath,amsthm,amsfonts,mathrsfs}
\usepackage[usenames,dvipsnames]{color}
\usepackage[pagebackref = true, colorlinks = true, linkcolor = blue, citecolor = Green]{hyperref}
\usepackage[alphabetic,lite]{amsrefs}
\usepackage{enumitem}
\usepackage{amscd}   
\usepackage[all, cmtip]{xy} 
\usepackage{xfrac}
\usepackage[T1]{fontenc}

\DeclareFontEncoding{OT2}{}{} 

\usepackage[all]{xy}
\usepackage{fullpage}

\usepackage{color} 

\def\act#1#2%
  {\mathop{}%
   \mathopen{\vphantom{#2}}^{#1}%
   \kern-3\scriptspace%
   #2}

\newcommand{\Q}{{\mathbb Q}}

\newcommand{\F}{{\mathbb F}}
\newcommand{\A}{{\mathbb A}}
\newcommand{\PP}{{\mathbb P}}

\newtheorem{Theorem}{Theorem}[section]

\newtheorem{Question}[Theorem]{Question}

\theoremstyle{definition}

\numberwithin{equation}{section}

\begin{document}
\title{Commitment Schemes and Diophantine Equations}

\author{Jos\'e Felipe Voloch}
\address{School of Mathematics and Statistics, University of Canterbury, Private Bag 4800, Christchurch 8140, New Zealand}
\email{felipe.voloch@canterbury.ac.nz}
\urladdr{http://www.math.canterbury.ac.nz/\~{}f.voloch}

\begin{abstract}
Motivated by questions in cryptography, we look for diophantine equations that are hard to solve but for which
determining the number of solutions is easy.
\end{abstract}

\maketitle

\section{Commitment Schemes}
\label{commit}

Solving a diophantine equation is typically hard but, given a point, it is typically easy to find a variety containing that
point. This is an example of a ``one-way function'' with potential applications to cryptography. Our current (lack of) 
knowledge suggests that such a function
is possibly quantum resistant and, therefore, cryptosystems based on these 
could be used for post-quantum cryptography \cite{BL}.

An encryption system based on this principle was proposed by Akiyama and Goto \cites{AG, AG2}, then broken by Ivanov and the author \cite{IV}. It was then fixed, broken again, fixed again,... Current status unclear.

The purpose of a commitment scheme is for a user to commit to a message without revealing it 
(e.g. vote, auction bid) by making public a value
obtained from the message in such a way that one can check, after the message is revealed, 
that the public value confirms the message.

Using such diophantine one-way functions for commitment schemes was proposed by Boneh and Corrigan-Gibbs \cite{BC}. 
They also suggested to work modulo an RSA modulus $N$. This could conceivably weaken the system. It will definitely no longer be quantum resistant.
Some partial attacks on this particular system are presented in \cite{ZW}.

Here is the general format of a diophantine commitment scheme.
Encode a message as point $P$ over some field $F$. Make public a variety $V/F$ with $P \in V$, 
with $V$ taken from some fixed family of varieties. To check the commitment, 
one verifies that $P$ satisfies the equations of $V$. We need the following conditions to be satisfied for
this to work:

\begin{itemize}
\item Given $P$, it is easy to construct $V$.
\item Given $V$, it is hard to find $V(F)$ (hence $P$).
\item Given $V$ (and perhaps $P$), it is easy to verify that $\# V(F) = 1$.
\end{itemize}

The last condition is important to prevent cheating. It proves that $P$ was indeed the committed message. In general, a
commitment scheme consists of two algorithms \textsf{Commit(m,r), Reveal(m,r,c)}. The first takes as input a message 
\textsf{m} and a random string \textsf{r} to produce an output \textsf{c}, which is then made public. The second takes as
input \textsf{m,r,c} as before and outputs \textsf{yes} or \textsf{no}, depending on whether \textsf{c} is the correct output of \textsf{Commit(m,r)}. The randomness is needed, e.g., if the list of possible messages is small enough that it can be brute
force searched.
See \cite[Section 4.1]{BC} for a more precise definition of a commitment scheme and some discussion.

These commitment schemes are similar in spirit to the class of
multivariate polynomial cryptosystems. In analogy to what is done
there, it is conceivable to have encryption by selecting a subset of 
varieties $V/F$ such that $V(F)$ can be easily found but that $V$ 
can be disguised as a general member of the collection of varieties. 
We do not address the interesting problem of doing this for schemes
we consider.

\section{Diophantine Equations}

Answering a question of Friedman, Poonen \cite{Poonen} proved:

\begin{Theorem} Assuming the Bombieri-Lang conjecture,
there exists $f(x,y) \in \Q[x,y]$
inducing an injection $\Q \times \Q \to \Q$.
\end{Theorem}

Boneh and Corrigan-Gibbs \cite{BC} then use the following construction from such a function.
For $P=(a,b)$, take $V: f(x,y)=f(a,b)$ to get a commitment scheme fitting the general setting of section \ref{commit}.  
Unfortunately, Poonen's proof, besides being conditional on a conjecture, is also non-constructive!

Zagier suggested $f(x,y)= x^7 + 3y^7$ as a polynomial defining an injective function. But we don't have a proof.
With exponent $13$ instead of $7$, the abcd conjecture implies that this function essentially injective. 

\begin{Question}
Is solving $x^7 + 3y^7=k$ over $\Q$ hard? 
\end{Question}

Pasten \cite{Pasten} proved that there exists an affine surface $S$ of the form $U \times U$ 
with $S(\Q)$ Zariski-dense in $S$
and a morphism $S \to \A^1$ inducing an injection $S(\Q) \to \Q$. But,  $S(\Q)$ is too sparse to be cryptographically useful.

Cornelissen \cite{Cornelissen}, using that the abcd conjecture is true for function fields of characteristic $0$,
noted that $x^m+ty^m$ is injective in $k(t), {\rm char} k  = 0$ for $m$ large.

\begin{Question}
Is solving $x^m + ty^m=k$ over $\Q(t)$ hard?
\end{Question}

My guess is that the answer is no.

He also noted that $x^p+ty^p$ is injective in $k(t), {\rm char} k  = p$.
But solving $x^p + ty^p=k$ is easy.

The following was noted in \cite{SV}, with the proof being an extension
of \cite{V} (see also \cite{Wang} for a related result without a hypothesis on the degree of the morphism):

\begin{Theorem}
\label{lem:S-unit eq}
Let $K$ be a function field of a curve $C$ of genus $g$ with field of constants $F$  of characteristic $p>0$ 
and let $S$ be a finite set of places of $K$.
If $u_1, \ldots, u_t$ are $S$-units of $K$, linearly independent over $F$,
such that the degree of the morphism $(u_1:\cdots:u_t): C \to \PP^{t-1}$ is less than $p$ and
satisfy
$$u_1 + \cdots  + u_t =1 $$
then 
$$
\max \{\deg u_i | i=1,\dots,t\} \le \frac{t(t-1)}{2} (2 g -2 +\# S)
$$
\end{Theorem}

The above result implies the injectivity
of $x^{13}+ty^{13}$ in the set of pairs of elements of $k(t)-k$ of degree
at most $p/13$ if $13 \nmid p(p-1)$.

This is enough for the application to commitment schemes by taking a sufficiently large finite field $k$ and considering the function $x^{13}+ty^{13}$ 
restricted to the above set where the function is injective.

But the function is not injective in the whole of $k(t)$. Indeed, if $x^{13}+ty^{13} = k, q=p^{12}$, then 
$$(x^q/k^{(q-1)/13})^{13} + t(t^{(q-1)/13}y^q/k^{(q-1)/13})^{13}=k$$

\section{Curves on surfaces}

The cryptosystem of Akiyama and Goto \cites{AG, AG2} actually uses
curves on surfaces over finite fields. We now consider the use 
of rational curves on surfaces in $\PP^3$ over a finite field for
commitment schemes.

We start with a rational curve $P$ parametrized by 
$(f_0:f_1:f_2:f_3)$
in $\PP^3$ over a finite field $\F_q$, where the $f_i$'s are 
polynomials of degree
at most $m$ (i.e a point in $\PP^3$ over $\F_q(t)$). 
Such a curve will include the message and randomness and our commitment will be
a smooth surface $S/\F_q$ of degree $d$ containing $P$. This is a bit different
from previous schemes as the surface is constant (i.e. independent of $t$). 
If $S$ is given by an homogeneous
equation $F=0$, the condition that $P \subset S$ is simply 
$F(f_0,f_1,f_2,f_3)=0$ which can be viewed as
a system of linear equation on the
coefficients of $F$, once the $f_i$ are given. There are ${d+3 \choose 3}$ coefficients
and $dm+1$ equations. One has solutions to the system as soon as there
are more coefficients than equations but these are not guaranteed to be
smooth. Poonen \cite{Poonen2} has proved that, for $d$ large, a positive
proportion of those solutions do indeed give smooth surfaces. One expects
in practice that, as long as the finite field is big enough, there will
be plenty of smooth surfaces.  

To guarantee uniqueness of the curve $P$ inside $S$, we prove the
following result.

\begin{Theorem}
Let $S/\F_q$ be a smooth surface in $\PP^3$ of degree $d>3$ with Picard 
number two. Then $S$ contains at most one smooth rational curve of degree
$m$, if $m < 2d(d-4)/(d-2)$.
\end{Theorem}

\begin{proof}
Let $H$ be a hyperplane section and $D_1,D_2$ two distinct 
smooth rational curves of degree $m$ contained in $S$. We compute the
determinant of
matrix of intersection pairings for $H,D_1,D_2$ and show it is non-zero, 
hence these curves are independent in the N\'eron-Severi group, contradicting
the hypothesis on the Picard number. 

Clearly, $H^2=d, HD_i = m, i=1,2$. The canonical class of $S$ is $(d-4)H$,
so the adjunction formula gives $D_i^2+(d-4)HD_i = -2$, hence
$D_i^2 = -(2+(d-4)m)$. Let $\delta = D_1D_2$. The 
determinant of matrix of intersection pairings is therefore
$$
\begin{vmatrix}
H^2 && HD_1 && HD_2 \\
D_1H && D_1^2 && D_1D_2 \\
D_2H && D_2D_1 && D_2^2 \\
\end{vmatrix}
=
\begin{vmatrix}
d && m && m \\
m && -(2+(d-4)m) && \delta \\
m && \delta && -(2+(d-4)m) \\
\end{vmatrix}
=
$$

$$-d\delta^2 +2m^2\delta +d(2+(d-4)m)^2+m^2(2+(d-4)m).$$
This vanishes precisely when $\delta = -(2+(d-4)m), 2m^2/d + (2+(d-4)m)$.
The first value is negative so cannot be $D_1D_2$ and the second
value is bigger than $m^2$ by our hypothesis but $D_1D_2 \le m^2$ by B\'ezout's
theorem so cannot be $D_1D_2$ either.

\end{proof}

To apply the theorem, we need to know that the Picard number of $S$ is
at most two. For a given
surface, this can be done using the algorithm of \cite{Costa}, for example.
This algorithm computes the $L$-function of $S$ and the 
Picard number of $S$ is the multiplicity of $q$ as a root of the 
$L$-function, conditional on the Tate conjecture. However, the surfaces
we construct will have Picard number at least two and a theorem of Tate
shows that the multiplicity of $q$ as a root of the
$L$-function is an upper bound for the Picard number. So, if this 
multiplicity is two, it is verified that the Picard number is two.
There is a parity condition coming from the functional equation for
$L$-functions which implies that this will not work if $d$ is odd.
It is reasonable to expect that a sizable proportion of such surfaces have 
Picard number two if $d$ is even, but this is not currently known and is worthy of further investigation.

In sum, our commitment scheme is as follows, with a finite field $\F_q$ and
integers $m,d$ selected a priori.

\begin{enumerate}
\item Encode a message as well as some randomness within $(f_0,f_1,f_2,f_3), f_i \in \F_q[t], \deg f_i \le m$.
\item Choose a random $F \in \F_q[x_0,x_1,x_2,x_3]$ homogeneous
of degree $d$ with $F(f_0,f_1,f_2,f_3)=0$.
\item Check whether the surface defined by $F=0$ is smooth and has Picard number two. If so, publish $F$ as the commitment. If not, pick a different $F$ in step (2).  
\end{enumerate}

For an explicit example, consider $m=3,d=6$. For a sextic surface to contain
a given twisted cubic, one needs to satisfy a system of $19$ equations in
$84$ variables and, hopefully, many of those will give rise to smooth surfaces
with Picard number two. The space of available messages depends on $16$ 
variables. 

One can also use $m=3,d=4$. The inequality in the theorem is not satisfied
but the second value for $\delta$ is $13/2$, which is not an integer so
cannot be $D_1D_2$ and the result holds. In this case, we have a system of
$13$ equations in $35$ variables for the coefficients of the surface and
again, the space of available messages depends on $16$
variables.

The expansion from $16$ variables to $84$ (or $35$) from the message to the
commitment is potentially wasteful and it is worth investigating whether a priori
setting many of these variables to zero will still allow enough variability so that step (3) 
above succeeds. Another issue worth studying is the choice of $q$. In some ways,
small $q$ is better for computations. But if a very small value of $q$, such as $q=2$ is
chosen, then $m=3$ is too small, as it allows brute force searching for the rational curve.

Given a surface, to find a 
rational curve inside it, one can either do a brute force search on the coefficients of the parametrization,
or set up a system of equations for these coefficients and try to solve it, e.g, using Gr\"obner bases. 
Neither option seem particularly efficient. Neither option also appears to be much improved by the use
of quantum computers.
There are general algorithms in the literature (e.g.
\cite{PTV}) that compute the N\'eron-Severi group of a variety but these
make no claim of practicality.

\section*{Acknowledgements}
This work was supported by MBIE. I would also like to thank Steven Galbraith for suggesting that I look into
commitment schemes and for helpful comments as well as Edgar Costa and Bjorn Poonen for suggestions.


\end{document}